\documentclass{report}

\usepackage[english]{babel}
\usepackage{enumitem}
\usepackage{manyfoot}
\usepackage[square,numbers]{natbib}
\DeclareNewFootnote{H}[Alph]

\usepackage{amsmath}
\usepackage{xcolor}
\usepackage{graphicx}
\usepackage{todonotes}
\usepackage{notation}

\usepackage{mathtools}
\mathtoolsset{showonlyrefs}

\newcommand{\dash}{\,---\,}

\title{Metric entropy of causal, discrete-time LTI systems}
\author{Clemens Hutter, Thomas Allard, and Helmut Bölcskei}

\begin{document}
\maketitle

\section{Introduction}
In \cite{hutter_metric_2022} it is shown that recurrent neural networks (RNNs) can learn\dash{}in a metric entropy optimal manner\dash{}discrete time, linear time-invariant (LTI) systems. This is effected by comparing the number of bits needed to encode the approximating RNN to the metric entropy of the class of LTI systems under consideration \cite{zamesConti, zamesDisc}.
The purpose of this note is to provide an elementary self-contained proof of the metric entropy results in \cite{zamesConti, zamesDisc}, in the process of which minor mathematical issues appearing in \cite{zamesConti, zamesDisc} are cleaned up. These corrections also lead to the correction of a constant in a result in \cite{hutter_metric_2022} (see Remark \ref{rem:correct_hutter}).

\paragraph*{Outline of the note.} Section \ref{sec:setting_and_statement} details the setup underlying the note and states the main result in Theorem \ref{thm:my_entropy}. The technical elements employed in the proof of Theorem \ref{thm:my_entropy} are presented in Section \ref{sec:covering_and_packing}.
For completeness, Appendix \ref{sec:metric_entropy_def} provides a brief introduction to the concept of metric entropy.

\paragraph{Notation.}
$\ind{\cdot}$ denotes the truth function which takes on the value $1$ if the statement inside $\{\cdot\}$ is true and equals $0$ otherwise.
For functions $f(\epsilon)$, $g(\epsilon)$, we use the notation $g(\epsilon) = o(f(\epsilon))$ to express that $\lim_{\epsilon\rightarrow 0 }\left|\frac{g(\epsilon)}{f(\epsilon)}\right| = 0$. $\log(\cdot)$ refers to the natural logarithm and $\log_2(\cdot)$ designates the logarithm to base 2.

\section{Setting and main statement}\label{sec:setting_and_statement}

The input-output relation of the LTI system $\sysOp$ is given by the 
convolution of the input signal $x[\cdot]$ with the system's impulse response $k[\cdot]$ according to
\begin{equation}\label{eq:sys_is_conv}
    (\sysOp x) [t] = \sum_{\tau=0}^{\infty} \impResp[\tau]x[t-\tau] =: (\impResp *x)[t], \qquad t\in\N_0,
\end{equation}
where we assume that $x[t]=0,\,k[t]=0$, both for $t < 0$, that is we consider one-sided input signals and causal systems. We occasionally use the notation $\impResp_\sysOp$ to designate the impulse response associated with the system $\sysOp$.
We shall frequently make use of the one-sided $\mathcal{Z}$-transform for $\ell_\infty$-signals defined as\footnote{
Note the positive exponents of $z$ in the definition of the $\mathcal{Z}$-transform. This convention is chosen to maintain consistency with Definition \ref{def:hardy} below adopted from \cite{zamesDisc}. 
} 
\begin{equation}\label{eq:z_def}
     X(z) := (\Z{x})(z) = \sum_{t=0}^{\infty} x[t]z^t,  \quad |z| < 1.
\end{equation}
Note that thanks to $x\in\ell_\infty$ the series \eqref{eq:z_def} converges absolutely for all $z\in\C$ with $|z|<1$.
Further, we will need the following norms. 
\begin{definition}[{\cite[Chapter~17]{Rudin1987}}] \label{def:hardy} 
      With $\setOfTransferFs := \{\Z{x} \mid x \in \ell_\infty, \, x[t] = 0,\text{ for } t<0 \}$,
      define for 
      $X \in \setOfTransferFs $ the Hardy norms \begin{align}
              \norm{X}_{\Hardy{2}} &:= \sqrt { \sup_{r\in(0,1)} \frac{1}{2\pi}\int_0^{2\pi} |X(re^{i\theta})|^2 d\theta, }\\
              \norm{X}_{\Hardy{\infty}} &:= \sup_{|z|<1} |X(z)|.
          \end{align}
      The corresponding Hardy spaces are given by $\Hardy{2}=\{X(\cdot)\,|\, X \in \setOfTransferFs, \; \norm{X}_{\Hardy{2}} < \infty  \}$ and $\Hardy{\infty}=\{X(\cdot)\,|\, X\in\setOfTransferFs, \; \norm{X}_{\Hardy{\infty}} < \infty\}$.
\end{definition}
Next, we note the well-known relation
\begin{equation}\label{eq:conv_is_multiplication}
    \left(\Z{\mathcal{L} x}\right)(z)=(\Z{k*x})(z) = K(z) \cdot X(z),
\end{equation}
where 
$K(z) := \left(\Z{\impResp}\right)(z)$ 
is commonly referred to as the system's transfer function. We now define the distance between LTI systems as follows. 

\begin{definition}\label{def:metric}
      For LTI systems $\sysOp$ and $\sysOp'$ with transfer functions $\transF(z)$ and $\transF'(z)$, respectively, both in $\Hardy{\infty}$, we define the metric
      \begin{align*}
            \rho(\sysOp, \sysOp') :=  \norm{\transF - \transF'}_{\Hardy{\infty}}. \label{eq:norm_def}
      \end{align*}      
\end{definition}
\begin{remark}
      By Theorem \ref{thm:norm_equivalence}, $\rho(\sysOp, \sysOp')$ can be expressed in the following forms
      \begin{align*}
            \rho(\sysOp, \sysOp') =&\;
            \norm{\transF - \transF'}_{\Hardy{\infty}}  \\
            =&\; \sup_{X\,\in\,\Hardy{2}} \frac{\norm{(\transF-\transF')X}_{\Hardy{2}}}{\norm{X}_{\Hardy{2}}}  \\
            =&\; \sup_{\norm{x}_{\ell^2} = 1} \norm{(k-k')* x}_{\ell^2}\\
            =&\;  \sup_{\norm{x}_{\ell^2} = 1} \norm{\sysOp x - \sysOp' x}_{\ell^2}.
      \end{align*}
\end{remark}

We are now ready to state the main result.
\begin{theorem}\label{thm:my_entropy}
      Let $\mC{}, \eC{} > 0$ and consider the set 
      \begin{equation}
          \mathcal{C}(\mC{}, \eC{})  = \setToCoverDef{}. \label{eq:setC}
      \end{equation}
      The metric entropy of $\mathcal{C}(\mC{}, \eC{})$ with respect to the metric
      \begin{equation*}
          \rho(\sysOp, \sysOp')  = \norm{\transF - \transF'}_{\Hardy{\infty}} 
          = \sup_{\norm{x}_{\ell^2} = 1} \norm{(k-k')* x}_{\ell^2}
      \end{equation*}
      satisfies
      \begin{equation}\label{eq:our_rate}
          \mEnt(\epsilon; \mathcal{C}(\mC{}, \eC{}), \rho) \thicksim \frac{\gamma}{2\eC{}}\left(\log\left(\frac{\mC{}}{\epsilon}\right)\right)^2,
      \end{equation}
      where $\gamma := \log_2(e)$ and $f(\epsilon)\thicksim g(\epsilon)$ stands for $\lim_{\epsilon \rightarrow 0} \left|\frac{f(\epsilon)}{g(\epsilon)}\right| = 1$.
      \begin{proof}
            Fix $\epsilon >0$. 
            According to Lemma \ref{lem:packing} below, there exists a $(2\epsilon)$-packing of $\mathcal{C}(\mC{}, \eC{})$ with $M_{2\epsilon}$ elements, where
            \begin{equation*}
                  \logII{M_{2\epsilon}} \geq \frac{\gamma}{2\eC{}}\left(\log\left(\frac{\mC{}}{\epsilon}\right)\right)^2 - o\left(\left(\log\left(\frac{\mC{}}{\epsilon}\right)\right)^2\right)
                  . 
            \end{equation*}
            Further, by Lemma \ref{lem:covering} below, there exists an $\epsilon$-covering of $\mathcal{C}(\mC{}, \eC{})$ with $N_\epsilon$ elements, where 
            \begin{equation*}
                  \logII{ N_\epsilon } \leq \frac{\gamma}{2\eC{}}\left(\log\left(\frac{\mC{}}{\epsilon}\right)\right)^2 + o\left(\left(\log\left(\frac{\mC{}}{\epsilon}\right)\right)^2\right)
                  .
            \end{equation*}
            Using Lemma \ref{lem:covering-packing}, we can hence sandwich the metric entropy according to 
            \begin{align*}
                  \frac{\gamma}{2\eC{}}\left(\log\left(\frac{\mC{}}{\epsilon}\right)\right)^2 - o\left(\left(\log\left(\frac{\mC{}}{\epsilon}\right)\right)^2\right)
                  &\leq
                  \logII{ M_{2\epsilon} } \\
                  &\leq
                  \logII{ M(2\epsilon; \mathcal{C}(\mC{}, \eC{}), \rho) } \\
                  &\leq 
                  \logII{N(\epsilon; \mathcal{C}(\mC{}, \eC{}), \rho)} \\
                 \left(\right . &=  \left . \mEnt{}(\epsilon; \mathcal{C}(\mC{}, \eC{}), \rho) \; \right) \\
                  & \leq \logII {N_\epsilon} \\
                  & \leq \frac{\gamma}{2\eC{}}\left(\log\left(\frac{\mC{}}{\epsilon}\right)\right)^2 + o\left(\left(\log\left(\frac{\mC{}}{\epsilon}\right)\right)^2\right).
            \end{align*}
            Dividing by $\frac{\gamma}{2\eC{}}\left(\log\left(\frac{\mC{}}{\epsilon}\right)\right)^2$ and taking $\lim_{\epsilon\rightarrow 0}$, then yields
            \[
                  1 \leq \lim_{\epsilon\rightarrow 0} 
                  \frac{ \mEnt{}(\epsilon; \mathcal{C}(\mC{}, \eC{}), \rho)}
                  {\frac{\gamma}{2\eC{}}\left(\log\left(\frac{\mC{}}{\epsilon}\right)\right)^2} \leq 1,
            \]
            which establishes \eqref{eq:our_rate}.
      \end{proof}
\end{theorem}

\begin{remark}\label{rem:correct_hutter}
      It follows from Theorem \ref{thm:my_entropy} that the constant $\frac{1}{b}$ in the scaling result 
      $\mEnt(\epsilon; \mathcal{C}(\mC{}, \eC{}), \rho) \thicksim \frac{1}{\eC{}}\left(\log\left(\frac{\mC{}}{\epsilon}\right)\right)^2$
      specified in \cite[Theorem 3.2]{hutter_metric_2022} is incorrect and should be replaced by $\frac{\gamma}{2b}$, i.e.,
      $\mEnt(\epsilon; \mathcal{C}(\mC{}, \eC{}), \rho) \thicksim \frac{\gamma}{2\eC{}}\left(\log\left(\frac{\mC{}}{\epsilon}\right)\right)^2$. 
      Similarly, the number of required bits specified as
      $\frac{1}{\eC{}}\left(\log\left(\frac{\mC{}}{\epsilon}\right)\right)^2 + o\left(\left(\log\left(\frac{1}{\epsilon}\right)\right)^2\right)$ 
      in \cite[Theorem 4.1]{hutter_metric_2022} should be replaced by 
      $\frac{\gamma}{2\eC{}}\left(\log\left(\frac{\mC{}}{\epsilon}\right)\right)^2 + o\left(\left(\log\left(\frac{1}{\epsilon}\right)\right)^2\right)$. 
\end{remark}

\section{Covering and packing bounds}\label{sec:covering_and_packing}

\newcommand{\constA}{\ensuremath{C_1}}
\newcommand{\cordset}[1]{\mathcal{X}_{#1}}
\newcommand{\deltaAt}[1]{\delta_{#1}}
\newcommand{\cover}{\mathcal{C}}
\newcommand{\packKer}[1]{\widetilde{k}_{#1}}
\newcommand{\packSys}[1]{\widetilde{\mathcal{L}}_{#1}}
\newcommand{\countT}[1]{n_{#1}}

\begin{lemma}\label{lem:packing}
      Let $\mC{}, \eC{} > 0$ and consider the set 
      \[
            \mathcal{C}(\mC{}, \eC{})  = \setToCoverDef{}
      \] equipped with the metric
      \begin{equation*}
            \rho(\sysOp, \sysOp')  = \norm{\transF - \transF'}_{\Hardy{\infty}} .
      \end{equation*}
      For every $\epsilon \in (0, a)$, there exists a $(2\epsilon)$-packing of $\mathcal{C}(\mC{}, \eC{})$ 
      with $M_{2\epsilon}$\hspace{-0.05em} elements, where
      \begin{equation*}
            \logII{M_{2\epsilon}} \geq \frac{\gamma}{2\eC{}}\left(\log\left(\frac{\mC{}}{\epsilon}\right)\right)^2 - o\left(\left(\log\left(\frac{\mC{}}{\epsilon}\right)\right)^2\right),
      \end{equation*}
      with $\gamma := \log_2(e)$.
      \begin{proof}
            We explicitly construct a $(2\epsilon)$-packing as visualized in Figure \ref{fig:packing}.
            \begin{figure}[t]
                  \includegraphics[width=0.95\textwidth]{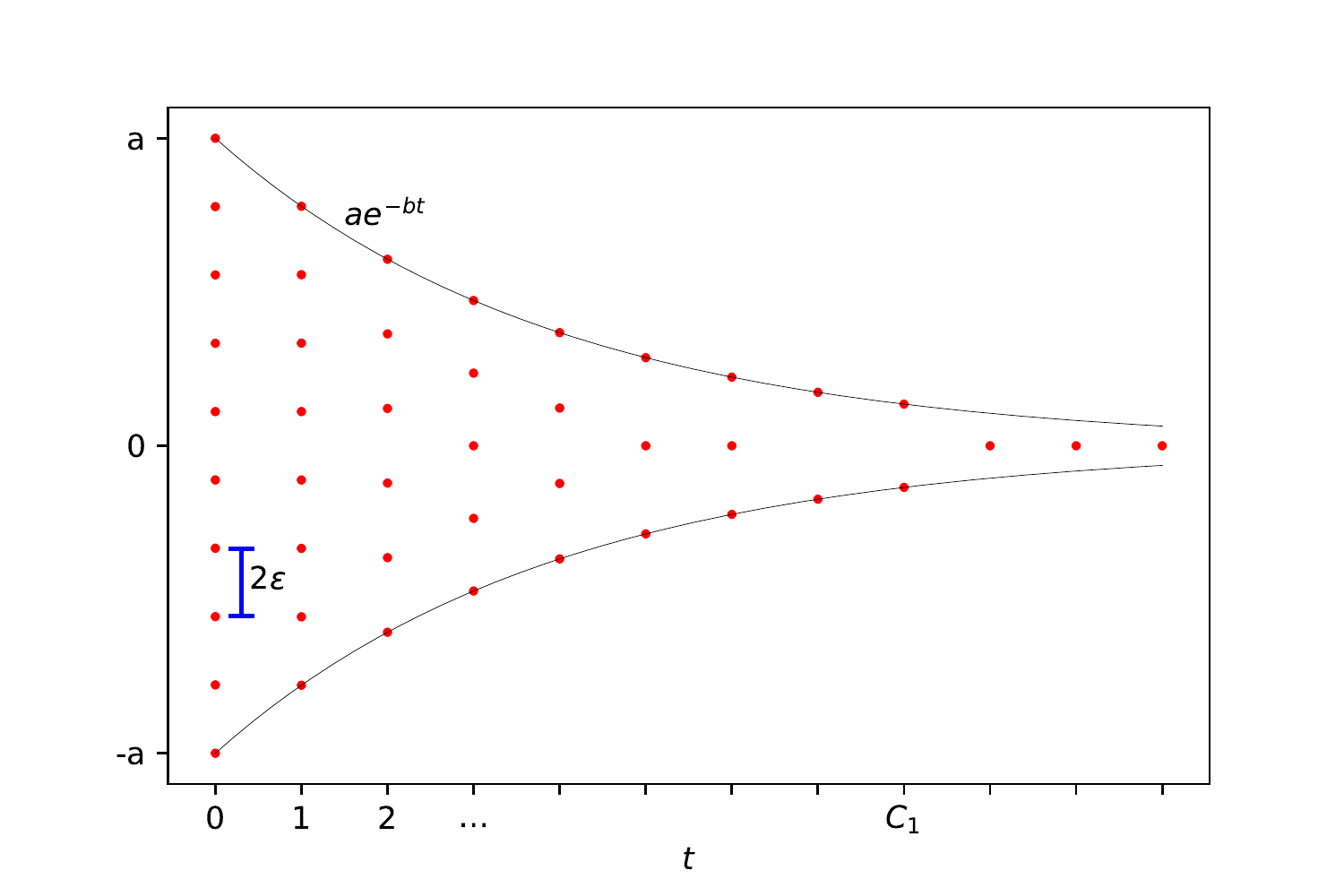}
                  \centering
                  \caption[Packing]{The packing is constructed by quantizing the impulse response at each time instant $t$ with quantization interval size of at least $2\epsilon$.}
                  \label{fig:packing}
            \end{figure}
            Define $\constA := \ceil*{\frac{1}{b} \log\left( \frac{a}{\epsilon}\right)} -1$. Now, for all  $  t \in \{0, \dots, \constA\}$, it holds that
            \[
            \begin{aligned}
                  & & t &\leq \constA < {\frac{1}{b} \log\left( \frac{a}{\epsilon}\right)} \\
                  &\Rightarrow \qquad &
                  b\,t &< \log\left( \frac{a}{\epsilon}\right) \\
                  &\Rightarrow \qquad &
                  -b\,t &> \log\left( \frac{\epsilon}{a}\right)  \\
                  &\Rightarrow \qquad &
                  e^{-bt} &> \frac{\epsilon}{a} \\
                  &\Rightarrow \qquad &
                  \frac{ae^{-bt}}{\epsilon} &> 1.
            \end{aligned}
            \]
            Next, for $t \in \{0, \dots, \constA\}$, we set  $ \countT{t} := \ceil*{\frac{2ae^{-bt}}{2\epsilon}}  - 1 \in \left[1, \frac{2ae^{-bt}}{2\epsilon}\right)  $ and $\deltaAt{t} := \frac{2ae^{-bt}}{ \countT{t}} > 2\epsilon$. 
            Now, we define the set
            \begin{equation*} \label{eq:packing}
                  \packingCollection := \left\{ \packSys{{\indI}_0, \dots, \indI_{\constA}}
                        \mid
                        \indI_\ell \in \{0, \dots,  \countT{\ell}\}, \text{ for } \ell \in \{0, \dots, \constA\} 
                  \right\},
            \end{equation*}
            where $\packSys{\indI_0, \dots, \indI_{\constA}}$ is the LTI system associated with the impulse response 
            \begin{equation*}
                  \packKer{\indI_0, \dots, \indI_{\constA}} [t] := 
                  \begin{cases} 
                        -ae^{-bt} + \indI_t \delta_t,  & 0 \leq  t  \leq \constA \\
                        0, & \text{otherwise} \\
                     \end{cases},
                     \label{eq:packing_def}
            \end{equation*}
            and we show that $\packingCollection$ constitutes a $(2\epsilon)$-packing of $\mathcal{C}(\mC{}, \eC{})$.
            First, we establish that $\packingCollection \subset  \mathcal{C}(\mC{}, \eC{})$ by verifying that  $\left| \packKer{\indI_0, \dots, \indI_{\constA}} [t] \right| \leq ae^{-bt}, \; \forall t\in\N_0$, holds for all $ \packKer{\indI_0, \dots, \indI_{\constA}} \in \packingCollection$.
            Indeed, for $t \in \{0, \dots, \constA\}$, we have
            \begin{equation*}
                  - ae^{-bt} \leq \packKer{\indI_0, \dots, \indI_{\constA}} [t] = -ae^{-bt} + \indI_t \delta_t \leq -ae^{-bt} +  \countT{t} \delta_t = -ae^{-bt} + 2ae^{-bt} = ae^{-bt},
            \end{equation*}
            and for $t > \constA$, 
            \begin{equation*}
                  \packKer{\indI_0, \dots, \indI_{\constA}} [t] = 0.  
            \end{equation*}

            \sloppy
            Next, we show that for distinct $\packSys{\indI_0, \dots, \indI_{\constA}},\, \packSys{j_0, \dots, j_{\constA}} \in \packingCollection$, i.e., there is at least one $\ell \in \{0, \dots, \constA\}$ such that $\indI_\ell \neq j_\ell$, it holds that \mbox{$\rho(\packSys{\indI_0, \dots, \indI_{\constA}}, \packSys{j_0, \dots, j_{\constA}})>2\epsilon$}.
            Indeed, for any such $\ell$, we have 
            \begin{align}
                  \rho(\packSys{\indI_0, \dots, \indI_{\constA}}, \packSys{j_0, \dots, j_{\constA}}) 
                  &=\sup_{\norm{x}_{\ell^2} = 1} \norm*{\left(\packKer{\indI_0, \dots, \indI_{\constA}}-\packKer{j_0, \dots, j_{\constA}}\right)* x}_{\ell^2} \label{teq:normeq} \\
                  & \geq \norm*{\packKer{\indI_0, \dots, \indI_{\constA}}-\packKer{j_0, \dots, j_{\constA}}}_{\ell^2}  \label{teq:chooseX} \\
                  & = \sqrt{ \sum_{t=0}^\infty \left(\packKer{\indI_0, \dots, \indI_{\constA}}[t]-\packKer{j_0, \dots, j_{\constA}}[t]\right)^2} \nonumber\\
                  & \geq \sqrt{ \left(\packKer{\indI_0, \dots, \indI_{\constA}}[\ell]-\packKer{j_0, \dots, j_{\constA}}[\ell]\right)^2} \nonumber\\
                  &= |-ae^{-b\ell} + \indI_\ell \delta_\ell + ae^{-b\ell} -j_\ell \delta_\ell | \nonumber\\
                  &= |\indI_\ell-j_\ell| \delta_\ell \nonumber\\
                  &> 2 \epsilon \label{teq:epsilon},
            \end{align}
            where in \eqref{teq:normeq} we used Theorem \ref{thm:norm_equivalence}, in \eqref{teq:chooseX} we inserted the particular choice $x[t] = \ind{t=0}[t]$ to lower-bound the $\sup$, and in \eqref{teq:epsilon} we used $\delta_\ell > 2\epsilon$.
            This establishes that $\packingCollection$ constitutes a $(2\epsilon)$-packing of $\mathcal{C}(\mC{}, \eC{})$ with respect to the metric $\rho(\cdot, \cdot)$ specified in Definition \ref{def:metric}. It remains to bound the cardinality of $\packingCollection$, which we denote by $M_{2\epsilon}$. Specifically, we have
            \begin{align}
                \logII{M_{2\epsilon}}  &= \log_2 \prod_{t=0}^{\constA{}} (1+  \countT{t})\nonumber\\
                &= \sum_{t=0}^{\constA} \log_2 \ceil*{\frac{a}{\epsilon}\,e^{-bt}}  \nonumber\\
                &\geq \sum_{t=0}^{\constA} \log_2 \left({\frac{a}{\epsilon}\,e^{-bt}}\right) \nonumber\\
                &= (\constA + 1) \log_2 \left(\frac{a}{\epsilon}\right) + \sum_{t=0}^{\constA} \log_2 \left(  e^{-bt} \right) \nonumber\\
                &= (\constA + 1) \log_2 \left(\frac{a}{\epsilon}\right) -b\;\log_2  (  e  ) \sum_{t=0}^{\constA} t  \, \nonumber\\
                &= (\constA + 1) \log_2 \left(\frac{a}{\epsilon}\right) -b\,\gamma\,\frac{\constA{}(\constA{}+1)}{2} \, \label{teq:log2e}\\
                &= \ceil*{\frac{1}{b} \log\left( \frac{a}{\epsilon}\right)} \log_2 \left(\frac{a}{\epsilon}\right) -b\,\gamma\,\frac{(
                      \ceil*{\frac{1}{b} \log\left( \frac{a}{\epsilon}\right)}-1)\ceil*{\frac{1}{b} \log\left( \frac{a}{\epsilon}\right)}}{2} \, \nonumber\\
                &= \ceil*{\frac{1}{b} \log\left( \frac{a}{\epsilon}\right)} \log_2 \left(\frac{a}{\epsilon}\right) 
                  - \frac{b\,\gamma}{2} \ceil*{\frac{1}{b} \log\left( \frac{a}{\epsilon}\right)}^2
                  + \frac{b\,\gamma}{2} \ceil*{\frac{1}{b} \log\left( \frac{a}{\epsilon}\right)} \nonumber\\
                &\geq 
                \frac{1}{b} \log\left( \frac{a}{\epsilon}\right) \log_2 \left(\frac{a}{\epsilon}\right) 
                - \frac{b\,\gamma}{2} \left( \frac{1}{b} \log\left( \frac{a}{\epsilon}\right) + 1 \right) ^2
                + \frac{b\,\gamma}{2} \frac{1}{b} \log\left( \frac{a}{\epsilon}\right) \nonumber\\
                &=
                \frac{1}{b} \log\left( \frac{a}{\epsilon}\right) \log_2 \left(\frac{a}{\epsilon}\right) 
                  - \frac{b\,\gamma}{2b^2} \left( \log\left( \frac{a}{\epsilon}\right) \right) ^2  
                  - \gamma \log\left( \frac{a}{\epsilon}\right)  
                  - \frac{b\,\gamma}{2}
                + \frac{\gamma}{2} \log\left( \frac{a}{\epsilon}\right) \nonumber\\
                &= \frac{\gamma}{b}\left(\log\left( \frac{a}{\epsilon}\right)\right)^2 
                - \frac{\gamma}{2b} \left( \log\left( \frac{a}{\epsilon}\right) \right) ^2  
                 - \frac{\gamma}{2} \log\left( \frac{a}{\epsilon}\right) 
                 -  \frac{b\,\gamma}{2}  \label{teq:log22logE}\\
                &= 
                 \frac{\gamma}{2b} \left( \log\left( \frac{a}{\epsilon}\right) \right) ^2  
                 - \frac{\gamma}{2} \log\left( \frac{a}{\epsilon}\right) 
                 -  \frac{b\,\gamma}{2}  \label{eq:packing_precise}\\
                &= \frac{\gamma}{2b} \left( \log\left( \frac{a}{\epsilon}\right) \right) ^2 - o\left(\left( \log\left( \frac{a}{\epsilon}\right) \right)^2\right), \nonumber
            \end{align}
            where 
            in \eqref{teq:log22logE} we used
            $\log_2(x) =\gamma \log(x)$. 
      \end{proof}
\end{lemma}

Before providing a covering for $\mathcal{C}(\mC{}, \eC{})$, we need the following auxiliary result.
\begin{lemma}\label{lem:hard_diff_l1_bound}
      Consider the LTI systems with impulse responses $k[\cdot]$ and $\widetilde{k}[\cdot]$ and corresponding transfer functions $K(z)$ and $\widetilde{K}(z)$, both in $\Hardy{\infty}$. We have
      \begin{equation*}
          \norm{K - \widetilde{K}}_{\Hardy{\infty}}  \leq \sum_{t=0}^{\infty} |k[t]-\widetilde{k}[t]|.
      \end{equation*}
      \begin{proof}
      The proof is by the following chain of relations 
      \begin{align*}
        \norm{K- \widetilde{K}}_\Hardy{\infty} &= 
        \sup_{|z|<1}\left| \sum_{t=0}^\infty k[t]z^{t} - \sum_{t=0}^\infty \widetilde{k}[t]z^{t} \right| \\
        & = \sup_{|z|<1}\left| \sum_{t=0}^\infty (k[t] - \widetilde{k}[t])z^{t} \right| \\
        &\leq \sup_{|z|<1}  \sum_{t=0}^\infty |k[t] - \widetilde{k}[t] | |z|^{t} \\
        &= \sum_{t=0}^\infty |k[t]- \widetilde{k}[t]|. &&\qedhere
       \end{align*}
      \end{proof}
  \end{lemma}

\newcommand{\constM}{C_2}
\newcommand{\coverKer}[1]{\packKer{#1}}
\newcommand{\coverSys}[1]{\packSys{#1}}
\newcommand{\gridMap}[1]{f_{#1}}
We are now ready to provide an upper bound on the covering number of $\mathcal{C}(\mC{}, \eC{})$.

\begin{lemma}\label{lem:covering}
      Let $\mC{}, \eC{} > 0$ and consider the set 
      \[
            \mathcal{C}(\mC{}, \eC{})  = \setToCoverDef{},
      \] equipped with the metric
      \begin{equation*}
            \rho(\sysOp, \sysOp')  = \norm{\transF - \transF'}_{\Hardy{\infty}}. 
      \end{equation*}
      For every $\epsilon \in (0, a) $, there exists an $\epsilon$-covering of $\mathcal{C}(\mC{}, \eC{})$ 
      with $N_{\epsilon}$ elements, where
      \begin{equation}
           \log_2 \left( N_{\epsilon} \right) \leq \frac{\gamma}{2\eC{}}\left(\log\left(\frac{\mC{}}{\epsilon}\right)\right)^2 + o\left(\left(\log\left(\frac{\mC{}}{\epsilon}\right)\right)^2\right),
      \end{equation}
      with $\gamma := \log_2(e)$.

      \begin{proof} 
            The proof is effected by explicit construction of an $\epsilon$-covering as visualized in Figure \ref{fig:covering}.
            \begin{figure}[t]
                  \includegraphics[width=0.95\textwidth]{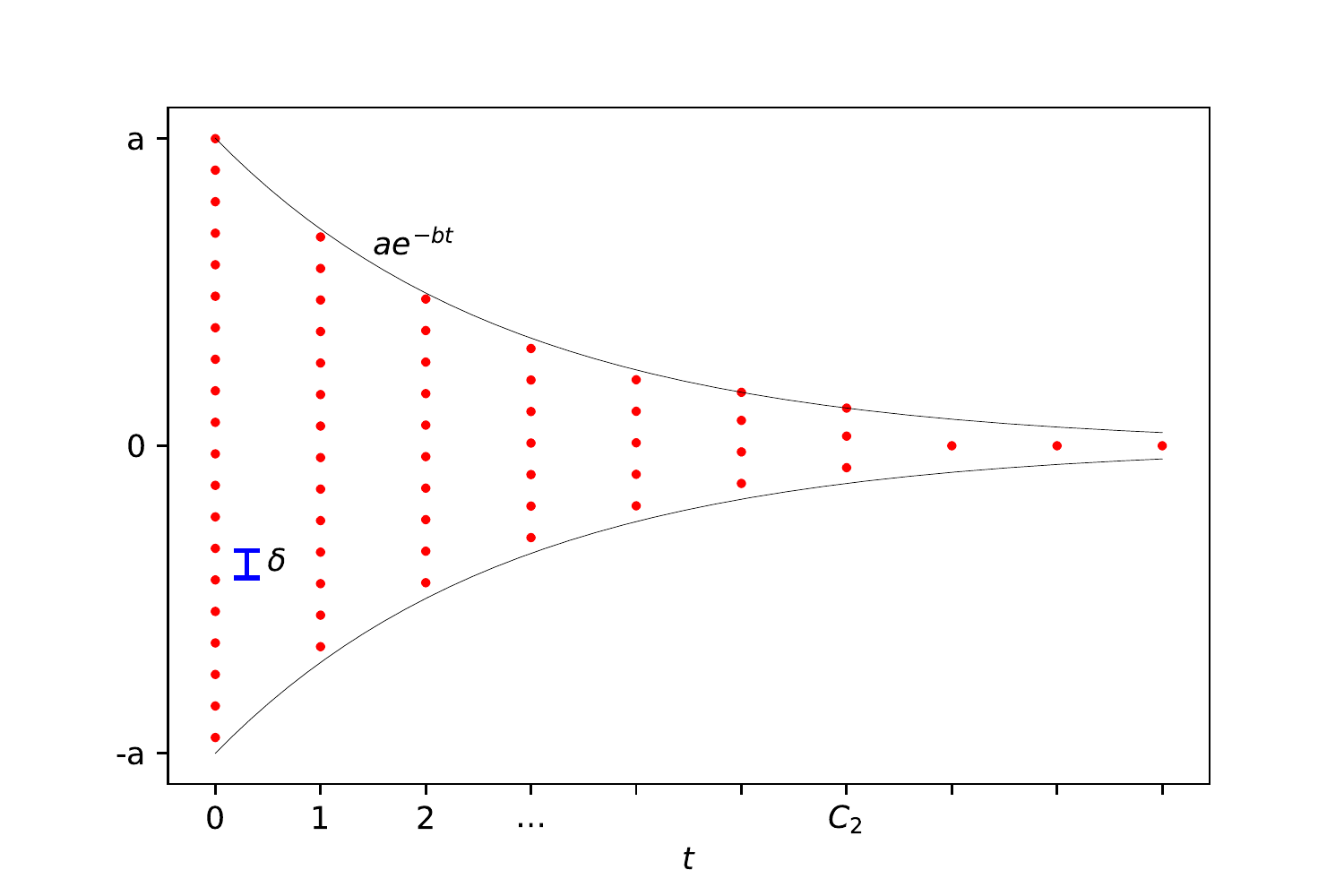}
                  \centering
                  \caption{
                  The covering is constructed by quantizing the impulse response at each time instant $t$ with quantization interval size $\delta$.
                  }
                  \label{fig:covering}
            \end{figure}
            We start by defining
            \begin{equation}
                  \constM:= \floor*{\frac{1}{\eC{}}\log \left(\frac{2a}{\epsilon\left(1-e^{-\eC}\right)} \right) }, 
                  \quad \delta := \frac{\epsilon}{\constM + 1}, 
                  \label{eq:define_cover_delta}
            \end{equation} 
            and 
            \[
                  \countT{t} := \ceil*{\frac{2ae^{-bt}}{\delta}}, \qquad \text{for } t \in \{0,\dots,\constM\}.
            \]
            As indicated in Figure \ref{fig:covering}, we quantize the impulse response at each time instant $t\in\{0,\dots, \constM\}$ 
            with quantization interval size $\delta$ using $n_t$ points. 
            To formalize this, we start by defining, for $t \in \{0, \dots, \constM\},$ the mappings
            \begin{align}
                  \begin{split}
                  \gridMap{t}:\; \{1, \dots, \countT{t}\} &\rightarrow [-ae^{-bt}, ae^{-bt}]\\
                  \label{eq:grid_map}
                  \gridMap{t}(i) &:= \min\left\{ -ae^{-bt} - \frac{\delta}{2} + i\delta,\, ae^{-bt} \right\} \,,
                  \end{split}
            \end{align}
            and show that, for all  $t \in \{0, \dots, \constM\}$, the following properties hold:
            \begin{enumerate}[label={\,(P\arabic*)}, ref={P\arabic*}]
                  \item $\gridMap{t}(1) \leq -ae^{-bt} + \frac{\delta}{2}$, \label{teq:p1}
                  \item $\gridMap{t}(\countT{t}) \geq ae^{-bt}  - \frac{\delta}{2}$,\label{teq:p2}
                  \item $\gridMap{t}(i+1) - \gridMap{t}(i) \leq \delta$, for $i\in \{1, \dots, \countT{t}-1\}$. \label{teq:p3}
            \end{enumerate}
            First, \eqref{teq:p1} follows by \[
                  \gridMap{t}(1) = \min\left\{ -ae^{-bt} - \frac{\delta}{2} + \delta,\, ae^{-bt} \right\} \leq 
                  -ae^{-bt} - \frac{\delta}{2} + \delta = -ae^{-bt} + \frac{\delta}{2}.
            \]
            To establish \eqref{teq:p2} we note that either
            \[
                  \gridMap{t}(\countT{t}) = ae^{-bt}
                  \qquad \text{ or } \qquad  
                  \gridMap{t}(\countT{t}) =  -ae^{-bt} - \frac{\delta}{2} + n_t\delta .
            \]
            In the former case we have $f_t(n_t) = ae^{-bt} \geq ae^{-bt} - \frac{\delta}{2}$, and in the latter, we obtain
            \[f_t(n_t) = -ae^{-bt} - \frac{\delta}{2} + n_t \delta \geq -ae^{-bt} - \frac{\delta}{2} + 2ae^{-bt} = ae^{-bt} - \frac{\delta}{2}.\]
            Finally, to prove \eqref{teq:p3}, we observe that, for all $ i \in \{1, \dots, \countT{t}-1\}$, it holds that
            \begin{align}
                   -ae^{-bt} - \frac{\delta}{2} + i\delta 
                 & \leq 
                   -ae^{-bt} - \frac{\delta}{2} + (\countT{t}-1)\delta \\
                &  \leq
                   -ae^{-bt} - \frac{\delta}{2} + \left( \frac{2ae^{-bt}}{\delta} \right)\delta\\
                   &= ae^{-bt} - \frac{\delta}{2} \leq ae^{-bt}.
            \end{align}
            Thus, $\gridMap{t}(i) = -ae^{-bt} - \frac{\delta}{2} + i\delta $, for all $ i \in \{1, \dots, \countT{t}-1\}$, and therefore
            \begin{align*}
                  \gridMap{t}(i+1) - \gridMap{t}(i) &= 
                  \min\left\{ -ae^{-bt} - \frac{\delta}{2} + (i+1)\delta,\, ae^{-bt} \right\} - \left( -ae^{-bt} - \frac{\delta}{2} + i\delta \right) \\
                  &\leq 
                  -ae^{-bt} - \frac{\delta}{2} + (i+1)\delta - \left( -ae^{-bt} - \frac{\delta}{2} + i\delta \right) = \delta,
                  \qquad  \text{ for } i \in \{1, \dots, \countT{t}-1\},
            \end{align*}
            which establishes \eqref{teq:p3}.

            Together, \eqref{teq:p1}-\eqref{teq:p3} imply that, for every $t\in\{0, \dots, \constM \}$ and $x \in [-ae^{-bt}, ae^{-bt}]$, there is an $i\in \{1, \dots, \countT{t}\}$ such that $|\gridMap{t}(i) - x | \leq \frac{\delta}{2}$. 
            Now, we define the set
            \begin{equation}
                  \coveringCollection{} := \left\{ \coverSys{{\indI}_0, \dots, \indI_{\constM}}
                        \mid
                        \indI_\ell \in \{1, \dots,  \countT{\ell}\}, \text{ for } \ell \in \{0, \dots, \constM\} 
                  \right\},
            \end{equation}
            where $\coverSys{\indI_0, \dots, \indI_{\constM}}$ is the LTI system associated with the impulse response 
            \begin{equation}
                  \coverKer{\indI_0, \dots, \indI_{\constM}} [t] := 
                  \begin{cases} 
                        \gridMap{t}(\indI_t),  & 0 \leq t \leq \constM \\
                        0, & \text{otherwise} \\
                     \end{cases} , \label{eq:cover_def}
            \end{equation}
            and we show that $\coveringCollection{}$ is, indeed, an $\epsilon$-covering for $\mathcal{C}(\mC{}, \eC{})$. Fix $\sysOp \in \mathcal{C}(\mC{}, \eC{})$ with corresponding impulse response $k[\cdot]$. As just established, for each $t\in \{0, \dots, \constM\}$, there is an $\indI_t \in \{1, \dots, \countT{t}\}$ such that $|k[t] - \gridMap{t}(\indI_t)| \leq \frac{\delta}{2}$. Hence, the corresponding $\coverKer{\indI_0, \dots, \indI_{\constM}} \in \coveringCollection{}$ satisfies 
            \begin{align}
                  \norm{K - \widetilde{K}_{\indI_0, \dots, \indI_{\constM}}}_{\Hardy{\infty}}  
                  & \leq \sum_{t=0}^{\infty} |k[t]-\widetilde{k}_{\indI_0, \dots, \indI_{\constM}}[t]| \label{teq:l1bound} \\
                    &= \sum_{t=0}^{\constM} |k[t]-\widetilde{k}_{\indI_0, \dots, \indI_{\constM}}[t]| + \sum_{t=\constM+1}^{\infty}|k[t]|  \\
                    &= \sum_{t=0}^{\constM} |k[t]-\gridMap{t}(\indI_t)| + \sum_{t=\constM + 1}^{\infty}|k[t]| \\
                    &\leq (\constM + 1) \frac{\delta}{2}   + \sum_{t=\constM+1}^{\infty}\mC{}e^{-\eC{} t} \label{eq:quant_iir_bounds}  \\
                    &= (\constM + 1) \frac{\delta}{2} + \mC{}\, \frac{e^{-\eC{}(\constM+1)}}{1-e^{-\eC}} \label{eq:both-bounds} \\
                    &\leq (\constM + 1) \frac{\delta }{2} + \mC{}\, \frac{e^{-\log\left(\frac{2\mC}{\epsilon(1-e^{-\eC})}\right)}}{1-e^{-\eC}} \label{eq:pluginM}\\
                    &= (\constM + 1) \frac{\delta }{2} + \mC{}\, \frac{\frac{\epsilon(1-e^{-\eC})}{2\mC}}{1-e^{-\eC}} \\
                    &=  \frac{\epsilon}{2} + \frac{\epsilon}{2} = \epsilon, \label{eq:quant_iir_final}
            \end{align}
            where in \eqref{teq:l1bound} we applied Lemma \ref{lem:hard_diff_l1_bound}, and in \eqref{eq:pluginM} we used $\constM \geq \frac{1}{\eC{}}\log \left(\frac{2a}{\epsilon\left(1-e^{-\eC}\right)} \right) - 1$.
            It remains to upper-bound $N_\epsilon$, the number of elements in $\coveringCollection{}$:
            \begin{align}
                  \log_2 \left( N_\epsilon \right) &= \log_2 \left( \prod_{t=0}^{\constM} n_t  \right) \nonumber \\
                        & = \sum_{t=0}^{\constM} \log_2 \left( n_t \right)  \nonumber \\
                        & = \sum_{t=0}^{\constM}\log_2 \left(  \ceil*{\frac{2ae^{-bt}}{\delta}} \right)  \nonumber\\
                        & \leq \sum_{t=0}^{\constM}\ceil*{\log_2 \left(  {\frac{2ae^{-bt}}{\delta}} \right)} \label{teq:ceilBoundLog} \\
                        & \leq \sum_{t=0}^{\constM} \left( {\log_2 \left(  {\frac{2ae^{-bt}}{\delta}} \right)} + 1 \right) \nonumber \\
                        &= -b\log_2(e)\sum_{t=0}^{\constM}  t  + (\constM+1)\left( \log_2\left(\frac{a}{\delta}\right) + 2\right)\nonumber\\
                        &= -\frac{\gamma b}{2}\constM(\constM+1) + (\constM+1)\left( \log_2\left(\frac{a}{\delta}\right) + 2\right)\nonumber\\
                        &= (\constM+1)\left(\gamma \log\left(\frac{a}{\delta}\right)  -\frac{\gamma b}{2}\,\constM + 2 \right)\label{teq:useLog22e}\\
                        &= (\constM+1)\left(\gamma \log\left(\frac{a(\constM+1)}{\epsilon}\right)  -\frac{\gamma b}{2}\,\constM + 2 \right)\nonumber\\
                        &= (\constM+1)\left(\gamma \log\left(\frac{a}{\epsilon}\right) + \gamma\log\left({\constM +1}\right)-\frac{\gamma b}{2}\,\constM + 2 \right) \nonumber,
            \end{align}
            where in \eqref{teq:ceilBoundLog} we used $\log_2\left(\ceil{x}\right) \leq \ceil*{\log_2(x)}$, $\forall x>0$, and in \eqref{teq:useLog22e} we employed $\log_2(x) =\gamma \log(x)$.
            Next, we note from the definition of $\constM$ that
            \begin{equation}
                 \frac{1}{\eC{}} \log \left(\frac{\mC{}}{\epsilon}\right)  + K_1(\eC{}) - 1 \leq \constM \leq \frac{1}{\eC{}} \log \left(\frac{\mC{}}{\epsilon} \right)  + K_1(\eC{}), \label{eq:upper-and-lower-on-m}
            \end{equation}
            with $K_1(\eC{}) := \frac{1}{\eC} \log \left(\frac{2}{1-e^{-\eC}} \right)$.
            Now we further upper-bound as follows:
            \begin{align}
                  &(\constM+1)\left(\gamma \log\left(\frac{a}{\epsilon}\right) + \gamma\log\left({\constM +1}\right)-\frac{\gamma b}{2}\,\constM + 2 \right)\\
                  &\leq (\constM+1) \left( 
                        \gamma \log\left(\frac{a}{\epsilon}\right)  
                        - \frac{\gamma b}{2}\left( \frac{1}{\eC{}} \log \left(\frac{\mC{}}{\epsilon}\right)  + K_1(\eC{}) - 1  \right) 
                        + \gamma\log\left({\constM +1}\right) + 2
                  \right) \\
                  &= (\constM+1) \left(
                        \frac{\gamma}{2}\log\left(\frac{a}{\epsilon}\right)  
                        + \gamma\log\left({\constM +1}\right) 
                        - \frac{\gamma b}{2} K_1(b) + \frac{\gamma b}{2}  + 2
                        \right) \\
                  &= (\constM+1) \left(
                        \frac{\gamma}{2}\log\left(\frac{a}{\epsilon}\right)  
                        + \gamma\log\left({\constM +1}\right) 
                        + K_2(b)
                        \right) \\
                  &\leq \left(\frac{1}{\eC{}} \log \left(\frac{\mC{}}{\epsilon} \right) + K_3(b) \right)
                         \left(
                              \frac{\gamma}{2}\log\left(\frac{a}{\epsilon}\right)  
                              + \gamma\log\left({\constM +1}\right) 
                              + K_2(b)
                        \right) \\
                  \begin{split}
                  &=
                  \frac{\gamma}{2b}  \left( \log \left(\frac{\mC{}}{\epsilon} \right)\right)^2 
                        + \frac{\gamma}{b} \log \left(\frac{\mC{}}{\epsilon} \right) \log(\constM+1)
                        + \frac{1}{b}  \log \left(\frac{\mC{}}{\epsilon} \right) K_2(b) \\ 
                  & \qquad 
                        + K_3(b ) \frac{\gamma}{2}\log\left(\frac{a}{\epsilon}\right)  
                        + K_3(b ) \gamma\log\left({\constM +1}\right) 
                        + K_3(b) K_2(b)
                  \end{split}
                        \\
                  \begin{split}
                  &=
                  \frac{\gamma}{2b}  \left( \log \left(\frac{\mC{}}{\epsilon} \right)\right)^2 
                        + \frac{\gamma}{b} \log \left(\frac{1}{\epsilon} \right) \log(\constM+1)
                        + K_4(b)\log \left(\frac{1}{\epsilon} \right) \label{eq:cover_precise}\\
                   &\qquad
                        + K_5(a,b)\log(\constM+1) + K_6(a,b) 
                  \end{split}
                        \\
                  &= 
                  \frac{\gamma}{2\eC{}}\left(\log\left(\frac{\mC{}}{\epsilon}\right)\right)^2 + o\left(\left(\log\left(\frac{\mC{}}{\epsilon}\right)\right)^2\right), 
            \end{align}
            where
            $K_2(b) := - \frac{\gamma b}{2} K_1(b) + \frac{\gamma b}{2}  + 2$, $K_3(b)  := K_1(b)+1$, $K_4(b) := \frac{1}{b} K_2(b) + \frac{\gamma}{2}K_3(b)$, {$K_5(a,b) := \frac{\gamma}{b}\log(a) + \gamma K_3(b)$}, and $K_6(a,b) := K_4(b) \log(a) + K_3(b)K_2(b)$. The last equality follows from $\log (\constM+1) = o(\log (\epsilon^{-1}))$.
      \end{proof}
\end{lemma}
\begin{corollary}
    Observing the precise nature of the lower and upper bounds in \eqref{eq:packing_precise} and \eqref{eq:cover_precise} respectively, we can also write 
    \begin{equation}\label{eq:our_rate_big_O}
          \mEnt(\epsilon; \mathcal{C}(\mC{}, \eC{}), \rho) =  \frac{\gamma}{2\eC{}}\left(\log\left(\frac{\mC{}}{\epsilon}\right)\right)^2 + \mathcal{O}\!\left( \log\left(\frac{1}{\epsilon}\right)\log \left( \log\left(\frac{1}{\epsilon} \right) \right)\right)
     \end{equation}
     instead of \eqref{eq:our_rate}.
\end{corollary}

\appendix
\section{Metric entropy}\label{sec:metric_entropy_def}

\begin{definition}[\cite{Wainwright2019}] \label{def:metric_entropy}
      Let $(\mathcal{X}, \rho)$ be  a metric space. An $\epsilon$-covering of a compact set $\mathcal{C} \subseteq \mathcal{X}$ with respect to the metric $\rho$ is a set of points $\{ x_1, \dots, x_N \} \subset \mathcal{C}$ such that for each $x \in \mathcal{C}$, there exists an $i \in [1, N]$ so that $\rho(x, x_i) \leq \epsilon$. The $\epsilon$-covering number $N(\epsilon; \mathcal{C}, \rho)$ is the cardinality of a smallest $\epsilon$-covering of $\mathcal{C}$ and $\mEnt(\epsilon; \mathcal{C}, \rho) := \log_2( N(\epsilon; \mathcal{C}, \rho))$ is the metric entropy of $\mathcal{C}$. 
\end{definition}

\begin{definition}
      Let $({\cal X},\rho)$ be a metric space.  An $\eps$-packing of a compact set ${\cal C} \subset {\cal X}$ with respect to the metric $\rho$ is a set $\{x_1,...,x_N\} \subset \cal C$ such that 
      $\rho(x_i,x_j) > \eps$, for all distinct $i,j$. The $\eps$-packing number $M(\eps; {\cal X},\rho)$ is the cardinality of the largest $\eps$-packing.
\end{definition}

\begin{lemma}\label{lem:covering-packing}
      Let $({\cal X},\rho)$ be a metric space and ${\cal C}$ a compact set in ${\cal X}$. For all $\eps > 0$, 
      it holds that
      $$
      M(2\eps; {\cal C},\rho) \le N(\eps; {\cal 
      C},\rho) \le M(\eps; {\cal C},\rho).
      $$
\begin{proof}
      First, choose a minimal $\eps$-covering and a maximal $2\eps$-packing of ${\cal C}$. Since no two centers of the $2\eps$-packing can lie in the same ball of the
      $\eps$-covering, it follows that $M(2\eps; {\cal C},\rho) \le N(\eps; {\cal C},\rho)$. To establish $N(\eps; {\cal C},\rho) \le M(\eps; {\cal C},\rho)$, we note that, given a maximum 
      $\epsilon$-packing with cardinality $M(\eps; {\cal C},\rho)$, for every $x\in {\cal C}$, we have the center of at least one of the balls in the packing within distance less than $\eps$. If this were not the case, we could
      add another ball to the packing thereby violating its maximality. This maximal packing hence also provides an $\eps$-covering 
      and hence 
      $N(\eps; {\cal C},\rho) \le M(\eps; {\cal C},\rho)$.
\end{proof}
\end{lemma}
\section{Norms}
The following discussion of the metric specified in Definition \ref{def:metric} largely follows \cite{hutter_metric_2022} and is reproduced here for completeness.

\newcommand{\countMeas}{\#}
\newcommand{\lebSig}{\mathcal{A}}
\newcommand{\countSpace}{(\N_0 \times \N_0)}
\newcommand{\countSig}{2^{(\N_0 \times \N_0)}}

\begin{theorem} \label{thm:z_isometry}
      Let $x\in \ell^2$ be a one-sided sequence, i.e., $x[t]=0$, for $t<0$. Then,
      we have
      \[
          \norm{X}_{\Hardy{2}} = \norm{x}_{\ell^2},
      \]
      where $X= \Z{x}$.
      \begin{proof}
          \begin{align} 
              \norm{X}_{\Hardy{2}}^2 &= \sup_{r\in(0,1)}\, \frac{1}{2\pi} \int_0^{2\pi} | X (re^{i\theta}) |^2 d\theta \\
               &= \sup_{r\in(0,1)}\,\frac{1}{2\pi} \int_0^{2\pi} \left| \sum_{t=0}^{\infty} x[t] (re^{i\theta})^t \right|^2 d\theta \\
               &= \sup_{r\in(0,1)}\,\frac{1}{2\pi} \int_0^{2\pi} \left( \sum_{t=0}^{\infty} x[t] (re^{i\theta})^t \right)  \overline{\left( \sum_{t'=0}^{\infty} x[t'] (re^{i\theta})^{t'} \right)} d\theta \\
               &= \sup_{r\in(0,1)}\,\frac{1}{2\pi} \int_0^{2\pi}  \sum_{t=0}^{\infty}\sum_{t'=0}^{\infty} x[t] r^t \; \overline{x[t']} r^{t'} \, e^{i\theta (t - t')}  d\theta \\
               &= \sup_{r\in(0,1)} \,  \sum_{t=0}^{\infty} \sum_{t'=0}^{\infty} x[t]  \overline{x[t']} \, r^{t+t'} \frac{1}{2\pi}\int_0^{2\pi} e^{i\theta (t - t')} d\theta \label{teq:use_fubini} \\
               &= \sup_{r\in(0,1)} \,  \sum_{t=0}^{\infty} \sum_{t'=0}^{\infty} x[t]  \overline{x[t']} \, r^{t+t'} \ind{t=t'} \\
               &= \sup_{r\in(0,1)} \,  \sum_{t=0}^{\infty} |x[t]|^2  \, r^{2t} \\
               &= \sum_{t=0}^{\infty} |x[t]|^2 = \norm{x}_{\ell^2}^2. &&\qedhere
          \end{align}
          In \eqref{teq:use_fubini}, we interchanged the order of integration and summation. 
          This step can be justified using the Fubini–Tonelli theorem as detailed in \cite[Section 2.6]{lapidoth_foundation_2017} and \cite[Exercise 10.J]{bartle_elements_1995}, since, $\forall r \in (0,1)$,
            \begin{align}
                  \frac{1}{2\pi} \int_0^{2\pi}  \sum_{t=0}^{\infty} \sum_{t'=0}^{\infty} \left|  x[t]\, r^t \;  \overline{x[t']}\,  r^{t'} e^{i\theta (t - t')} \right| d\theta 
                  &= \frac{1}{2\pi}  \int_0^{2\pi} \left( \sum_{t=0}^{\infty} |x[t]|r^t \right) \left( \sum_{t'=0}^{\infty}  |x[t']|r^{t'} \right) d\theta \\
                  &\leq \frac{1}{2\pi}  \int_0^{2\pi} \norm{x}_{\ell_\infty}^2 \left( \sum_{t=0}^{\infty} r^t \right)^2 d\theta \\
                  &= \norm{x}_{\ell_\infty}^2 \left(\frac{1}{1-r}\right)^2 \\
                  &\leq \norm{x}_{\ell_2}^2  \left(\frac{1}{1-r}\right)^2 < \infty.
            \end{align}
      \end{proof}
      \end{theorem}
      
      \begin{theorem}\label{thm:opNorm_eq_infty}
      For $K(\cdot)\in\Hardy{\infty}$, 
      it holds that
      \begin{equation}
           \hNorm{\transF} = \sup_{X\in \Hardy{2}} \frac{\hhNorm{\transF X}}{\hhNorm{X}}. \label{eq:Rochberg}
      \end{equation}
      \begin{proof}
      The proof essentially follows \cite{McCarthy03} with some details filled in and minor refinements.
      We start by noting that the RHS of \eqref{eq:Rochberg} is the operator norm $\opNorm{K}:=\sup_{X\in \Hardy{2}} \frac{\hhNorm{KX}}{\hhNorm{X}}$ of the multiplication operator $X(z) \rightarrow K(z) X(z)$ and first establish that $\opNorm{K} \leq \hNorm{\transF}$. Indeed, for every $X\in \Hardy{2}$,
      we have
              \begin{align*}
              \norm{K\, X}_\Hardy{2} &= \HardyTwo{|K(re^{i\theta})  X(re^{i\theta})|^2} \\
              & \leq \HardyTwo{|X(re^{i\theta})|^2  \left(\sup_{|z|<1} |K(z)| \right)^2} \\
              & = \norm{K}_\Hardy{\infty}  \HardyTwo{|X(re^{i\theta})|^2} \\
              &= \norm{K}_\Hardy{\infty} \norm{X}_\Hardy{2},
            \end{align*}
      which, upon division by $\norm{X}_\Hardy{2}$, establishes the desired upper bound $\opNorm{K} \leq \hNorm{\transF}$.
      
      To complete the proof, we show that $\opNorm{K} \geq \hNorm{\transF}$.
      Applying 
      \[\hhNorm{\transF X} \leq \opNorm{\transF} \hhNorm{X}
      \]
      repeatedly, we get, for every $n\in \N$,
      \begin{align}\label{teq:22}
          \hhNorm{K^n X} \leq \opNorm{K}^n \hhNorm{X}.
      \end{align}
      Without loss of generality, we can restrict ourselves to $\opNorm{\transF}=1$ as otherwise we can simply consider $\transF' := \transF / \opNorm{\transF}$. 
      Next, towards a contradiction, assume that $\opNorm{\transF} < \hNorm{\transF}$, which, thanks to $\opNorm{\transF}=1$, results in
      $1< \hNorm{\transF} = \sup_{r<1, \, 0\, \leq\, \theta < 2\pi} |\transF(re^{i\theta})|$. 
      As $\transF(\cdot) \in \Hardy{\infty }$ by assumption, it follows that $\transF(z)$ is analytic and thus continuous inside the unit disk. 
      Hence, there exist $r' \in (0, 1), \epsilon>0$ and an interval $[\underline{\theta}, \overline{\theta})\,\in\,[0,2\pi)$
      with $ \overline{\theta} - \underline{\theta}= \delta >0$ such that
      \begin{equation}\label{teq:a24}
      |\transF(r'e^{i\theta'})| > 1 + \epsilon, \quad \forall \theta' \in [\underline{\theta}, \overline{\theta}).
      \end{equation}
      Now we take $X(z)=1=\Z{\ind{t=0}[t]}$ 
      which clearly satisfies $\hhNorm{X}=1$. Inserting this into \eqref{teq:22}, we obtain
      \begin{equation*}
           \hhNorm{\transF^n X}^2 \leq \opNorm{\transF}^{2n}\hhNorm{X}^2 = 1.
      \end{equation*}
      This finalizes the proof by leading to the following contradiction
      \begin{align}
          1 &\geq \hhNorm{\transF^n X}^2 \nonumber \\
          &= \sup_{0<r<1} \frac{1}{2\pi} \int_0^{2\pi} |\transF(re^{i\theta})|^{2n} d\theta \nonumber\\
          & \geq \frac{1}{2\pi} \int_0^{2\pi} |\transF(r'e^{i\theta})|^{2n} d\theta \nonumber \\
          & \geq \frac{1}{2\pi} \int_0^{2\pi} ((1+\epsilon)\ind{\theta \in [\underline{\theta}, \overline{\theta})})^{2n} d\theta \label{teq:a27}\\
          &= \frac{\delta}{2\pi} (1+\epsilon)^{2n} \xrightarrow[n\rightarrow \infty]{} \infty, \nonumber
      \end{align}
      where in
      \eqref{teq:a27} we used \eqref{teq:a24}.
      \end{proof}
      \end{theorem}

\begin{theorem}\label{thm:norm_equivalence}
      Let $\sysOp$ and $\sysOp'$ be LTI systems with corresponding transfer functions $K(z)$ and $K'(z)$, both in $\Hardy{\infty}$. 
      It holds that
      \begin{align}
            \norm{\transF - \transF'}_{\Hardy{\infty}}  
            &= \sup_{X\,\in\,\Hardy{2}} \frac{\norm{(\transF-\transF')X}_{\Hardy{2}}}{\norm{X}_{\Hardy{2}}}  \label{teq:aa1}\\
            &= \sup_{\norm{x}_{\ell^2} = 1} \norm{(k-k')* x}_{\ell^2}\label{teq:aa2}\\
            &=  \sup_{\norm{x}_{\ell^2} = 1} \norm{\sysOp x - \sysOp' x}_{\ell^2}\label{teq:aa3}.
      \end{align}
      \begin{proof}
            Equation \eqref{teq:aa1}  follows from Theorem \ref{thm:opNorm_eq_infty} upon noting that $K-K'\,\in\,\Hardy{\infty}$ by application of the triangle inequality. Next, \eqref{teq:aa2} is established through 
            \begin{align}
            \sup_{X\,\in\,\Hardy{2}} \frac{\norm{(\transF-\transF')X}_{\Hardy{2}}}{\norm{X}_{\Hardy{2}}}
            &= 
            \sup_{x\in\ell_2} \frac{\norm{(k-k')* x}_{\ell^2} }{\norm{x}_{\ell_2}}\\
            &= \sup_{x\in\ell_2} {\norm*{(k-k')* \frac{x}{\norm{x}_{\ell_2}}}_{\ell^2} }\\
            &= \sup_{\norm{x}_{\ell_2} = 1} \norm{(k-k')* x}_{\ell^2},
            \end{align}
            where we used Theorem \ref{thm:z_isometry}, \eqref{eq:conv_is_multiplication}, and the fact that convolution is linear. Finally, \eqref{teq:aa3} follows directly from \eqref{eq:sys_is_conv}.
      \end{proof}
\end{theorem}

\bibliographystyle{customels}

\bibliography{main}

\begin{thebibliography}{8}
\expandafter\ifx\csname natexlab\endcsname\relax\def\natexlab#1{#1}\fi
\providecommand{\url}[1]{\texttt{#1}}
\providecommand{\href}[2]{#2}
\providecommand{\path}[1]{#1}
\providecommand{\DOIprefix}{doi:}
\providecommand{\ArXivprefix}{arXiv:}
\providecommand{\URLprefix}{URL: }
\providecommand{\Pubmedprefix}{pmid:}
\providecommand{\doi}[1]{\href{https://doi.org/#1}{\path{#1}}}
\providecommand{\Pubmed}[1]{\href{pmid:#1}{\path{#1}}}
\providecommand{\bibinfo}[2]{#2}
\ifx\xfnm\relax \def\xfnm[#1]{\unskip,\space#1}\fi
\bibitem[{Hutter et~al.(2022)Hutter, Gül, and Bölcskei}]{hutter_metric_2022}
\bibinfo{author}{C.~Hutter}, \bibinfo{author}{R.~Gül}, and
  \bibinfo{author}{H.~Bölcskei},
\newblock \bibinfo{title}{Metric entropy limits on recurrent neural network
  learning of linear dynamical systems},
\newblock \bibinfo{journal}{Applied and Computational Harmonic Analysis}
  \bibinfo{volume}{59} (\bibinfo{year}{2022}) \bibinfo{pages}{198--223}.
  \DOIprefix\doi{10.1016/j.acha.2021.12.004}.
\bibitem[{Zames(1979)}]{zamesConti}
\bibinfo{author}{G.~Zames},
\newblock \bibinfo{title}{{On the metric complexity of causal linear systems:
  $\epsilon$-entropy and $\epsilon$-dimension for continuous time}},
\newblock \bibinfo{journal}{IEEE Transactions on Automatic Control}
  \bibinfo{volume}{24} (\bibinfo{year}{1979}) \bibinfo{pages}{222--230}.
  \DOIprefix\doi{10.1109/TAC.1979.1101976}.
\bibitem[{Zames and Owen(1993)}]{zamesDisc}
\bibinfo{author}{G.~Zames} and \bibinfo{author}{J.~G. Owen},
\newblock \bibinfo{title}{{A note on metric dimension and feedback in discrete
  time}},
\newblock \bibinfo{journal}{IEEE Transactions on Automatic Control}
  \bibinfo{volume}{38} (\bibinfo{year}{1993}) \bibinfo{pages}{664--667}.
  \DOIprefix\doi{10.1109/9.250545}.
\bibitem[{Rudin(1987)}]{Rudin1987}
\bibinfo{author}{W.~Rudin}, \bibinfo{title}{Real and Complex Analysis},
  \bibinfo{edition}{3rd} ed., \bibinfo{publisher}{McGraw-Hill},
  \bibinfo{year}{1987}.
\bibitem[{Wainwright(2019)}]{Wainwright2019}
\bibinfo{author}{M.~J. Wainwright}, \bibinfo{title}{{High-Dimensional
  Statistics}}, \bibinfo{publisher}{Cambridge University Press},
  \bibinfo{year}{2019}. \DOIprefix\doi{10.1017/9781108627771}.
\bibitem[{Lapidoth(2017)}]{lapidoth_foundation_2017}
\bibinfo{author}{A.~Lapidoth}, \bibinfo{title}{A {Foundation} in {Digital}
  {Communication}}, \bibinfo{edition}{2} ed., \bibinfo{publisher}{Cambridge
  University Press}, \bibinfo{address}{Cambridge}, \bibinfo{year}{2017}.
  \URLprefix
  \url{https://www.cambridge.org/core/books/foundation-in-digital-communication/05F46005A017815C49810D51CFDB9B8F}.
  \DOIprefix\doi{10.1017/9781316822708}.
\bibitem[{Bartle(1995)}]{bartle_elements_1995}
\bibinfo{author}{R.~G. Bartle}, \bibinfo{title}{The elements of integration and
  {Lebesgue} measure}, \bibinfo{publisher}{Wiley}, \bibinfo{address}{New York},
  \bibinfo{year}{1995}.
\bibitem[{McCarthy(2003)}]{McCarthy03}
\bibinfo{author}{J.~E. McCarthy},
\newblock \bibinfo{title}{{Pick's Theorem-What's the big deal?}},
\newblock \bibinfo{journal}{The American Mathematical Monthly}
  \bibinfo{volume}{110} (\bibinfo{year}{2003}) \bibinfo{pages}{36--45}.
  \DOIprefix\doi{10.2307/3072342}.

\end{thebibliography}

\end{document}